\documentclass[12pt]{article}

\usepackage{amsmath}
\usepackage{amsfonts}
\usepackage[T1]{fontenc} 

\newtheorem{theo}{Theorem}[section]
\newtheorem{prop}{Proposition}[section]
\newtheorem{corol}{Corollary}[section]

\newtheorem{lemme}{Lemma}[section]
\newtheorem{remark}{Remark}[section]

\newtheorem{example}{Example}[section]

\newcommand{\CQFD}{\hfill $\square$}

\newcommand{\cB}{{\cal B}}

\newcommand{\cF}{{\cal F}}

\newcommand{\cK}{{\cal K}}
\newcommand{\cL}{{\cal L}}
\newcommand{\cM}{{\cal M}}

\newfont{\msbm}{msbm10 scaled\magstep1}%blackboardbold
\newfont{\msbms}{msbm7 scaled\magstep1} %blackboardbold   subscript

\newcommand{\bbC}{\mbox{$\mbox{\msbm C}$}}

\newcommand{\bbE}{\mbox{$\mbox{\msbm E}$}}

\newcommand{\bbN}{\mbox{$\mbox{\msbm N}$}}

\newcommand{\bbP}{\mbox{$\mbox{\msbm P}$}}

\newcommand{\bbR}{\mbox{$\mbox{\msbm R}$}}

\newcommand{\bbsR}{\mbox{$\mbox{\msbms R}$}}

\def\E{\mathop{\hbox{\rm I\kern-0.20em E}}\nolimits}

\def\og{\leavevmode\raise.3ex
     \hbox{$\scriptscriptstyle\langle\!\langle$~}}
\def\fg{\leavevmode\raise.3ex
     \hbox{~$\!\scriptscriptstyle\,\rangle\!\rangle$}~}
     
\begin{document}

\title{Extremal shot noises, heavy tails and max-stable random fields.}
\author{ 
Cl\'ement Dombry\footnote{Laboratoire LMA, Universit\'e de Poitiers, 
 Téléport 2, BP 30179, F-86962 Futuroscope-Chasseneuil cedex, France. 
 Email: clement.dombry@math.univ-poitiers.fr}
 \protect\hspace{1cm}}
%\date{}
\maketitle

\begin{abstract}
{We consider the extremal shot noise defined by
$$M(y)=\sup\{mh(y-x);(x,m)\in\Phi\},$$
where $\Phi$ is a Poisson point process on $\bbR^d\times (0,+\infty)$ with intensity $\lambda dxG(dm)$ and $h:\bbR^d\to [0,+\infty]$ is a measurable function. Extremal shot noises naturally appear in extreme value theory as a model for spatial extremes and serve as basic models for annual maxima of rainfall or for coverage field in telecommunications. In this work, we examine their properties such as boundedness, regularity and ergodicity. Connections with max-stable random fields are established: we prove a limit theorem when the distribution $G$ is heavy-tailed and the intensity of points $\lambda$ goes to infinity. We use a point process approach strongly connected to the Peak Over Threshold method used in extreme value theory. Properties of the limit max-stable random fields are also investigated.}
\end{abstract}


\section{Introduction}
In this work we study extremal shot noises and their properties, in connection with extreme value theory. 
Extremal shot noises are flexible models that arise naturally when dealing with extremal events in a spatial setting. 
They are defined by 
\begin{equation}\label{eq:defSN}
M(y)=\sup \left\{mh(y-x);\ (x,m)\in\Phi\right\},\quad y\in\bbR^d,
\end{equation}
where $\Phi$ is a Poisson point process on $\bbR^d\times (0,+\infty)$ with intensity $\lambda dxG(dm)$. Here, $\lambda>0$ is the intensity of points, $G$ denotes a measure on $(0,+\infty)$ and $h:\bbR^d\to [0,+\infty]$ is a measurable function referred to as the shape function. The model can easily be extended to random shape functions (some examples are given below) but we will limit our study to deterministic shape functions.

Let us mention two examples where such random fields naturally arise. The first one is the analysis of annual maxima of daily spatial rainfall. Modeling rainfall is a very complex task and a vast amount of literature on this topic exists both in statistics and applied sciences, see for example \cite{BN} and the references therein. Convective precipitation has usually a local area of high intensity and minor to no rainfall elsewhere, with a superposition of many storm events. To deal with this type of precipitation, Smith \cite{Sm} proposed the so-called storm process given by equation \eqref{eq:defSN} for a Poisson point process $\Phi$ on $\bbR^d\times (0,+\infty)$ with intensity measure $dxm^{-2}dm$. In this context, a point $(x,m)$ of the point process represents a storm event located at $x$ with intensity $m$. The function $h$ on $\bbR^d$ is  non-negative with $\int_{\bbsR^d} h(x)dx=1$ and gives the typical shape of a storm event. 
The process $M$ is then a stationary max-stable spatial process with unit Fréchet margins. In \cite{Sm}, Smith proposed the multivariate Gaussian density with covariance matrix $\Sigma$ as a typical choice for the shape function $h$.
To illustrate the flexibility of such models, consider that we wish to assign to each storm event a spatial extension given by a radius parameter $r>0$ chosen at random with distribution $F(dr)$. This yields the model
$$M(y)=\sup_{(x,m,r)\in\Phi} mh\left(\frac{y-x}{r}\right),\quad y\in\bbR^d,$$
where $\Phi$ is a Poisson point process on $\bbR^d\times (0,+\infty)\times (0,+\infty)$ with intensity $dxG(dm)F(dr)$. This is an instance of extremal shot noise with random shape function $h_r(\cdot)=h(r^{-1}\cdot)$. We could also wish to introduce the temporal dimension and assign to each storm event the time $u>0$ when it occurs. Considering the extremal process at point $y$ up to time $t$, we define 
$$M(y,t)=\sup\left\{ mh(y-x)1_{\{0\leq u\leq t\}}; (x,m,u)\in\Phi\right\},\quad (y,t)\in\bbR^d\times [0,+\infty),$$
where $\Phi$ a Poisson point process on $\bbR^d\times (0,+\infty)\times (0,+\infty)$ with intensity $dxG(dm)dt$. In \cite{Schla}, Schlater considered the case when the shape function $h$ is given by a random process, as well as another class of max-stable process more suitable to deal with cyclonic precipitation with variable rainfall all over the region.

Extremal shot noises also arise in the domain of telecommunications. In this context, the Poisson point process $\Phi$ stands for a set of transmitters in the Euclidean space. A point $(x,m)$ is seen as a transmitter located at position $x$ and with power $m$. The function $h$ is the so-called attenuation function such that $mh(y-x)$ stands for the signal power received at $y$ from the transmitter $(x,m)$. In this context, a typical choice for $h$ is the omni-directional path-loss function defined by
$$ h(u)=(A \max(r_0,|u|))^{-\beta}\quad \mbox{or}\quad  h(u)=(1+A|u|)^{-\beta}$$
for some $A>0$, $r_0>0$ and $\beta>d$, which is the so-called path-loss exponent. In this setting, the extremal shot noise $M(y)$ given by equation $\eqref{eq:defSN}$ stands for the maximal power transmission at location $y$. Note that this scenario is isotropic, i.e. all the antennas are omni-directional. A more realistic scenario with directional antennas could be described as follows: for simplicity, we consider the planar case $d=2$; the antenna azimuth is denoted by $\theta\in [0,2\pi)$ and is considered as an additional mark, so that each transmitter is now represented by a point $(x,m,\theta)$ of a Poisson point process on $\bbR^2\times(0,+\infty)\times [0,2\pi)$; the power received at point $y$ from transmitter $(x,m,\theta)$ is then given by $mh_\theta(y-x)$  with $h_\theta$ given for example by
$$h_\theta(u)=m\alpha^2(\theta-\arg(u))(1+A|u|)^{-\beta}$$
with $\alpha^2:[0,2\pi)\to \bbR^+$ the radiation pattern of the antennas (see \cite{BB}). Since the shape function $h_\theta$ depends on the random mark $\theta$, it can be seen as a random shape function.

The extremal shot noise model defined here is closely related to mixed maxima moving process \cite{Schla, ZS04, S08}; in these papers however, only max-stable random fields are investigated, corresponding to suitable choice of the weight distribution $G$. On the contrary in the present work, we consider a general distribution $G$: the max-stable property is lost, but the interesting property of max-infinitely divisibility remains true (see \cite{EHV90}). Interesting links with stochastic geometry and union shot noise of random closed sets are established in Heinrich and Molchanov  \cite{HM}; in \cite{RR91}, Resnick and Roy provide applications of the theory of random upper semi-continuous functions and max-stable processes to the continuous choice problem. 

Deep connections exist between extreme value theory and regular variations; they are now well known in the multivariate case \cite{BGT, Res}, whereas many recent works focus on the functional case \cite{dH, HL05, HL07}. For example, de Haan and Lin \cite{dH} characterized the domain of attraction of continuous max-stable processes on $[0,1]$; Davis and  Mikosch \cite{DM} considered the notion of regularly varying process in connection with moving average models and space-time max-stable random fields;  results on extremes of moving average driven by a general regularly varying L\'evy process are obtained by Fasen \cite{F05}; Kabluchko, Schlather and de Haan \cite{KSH}  studied max-stable random fields arising as extremes of Gaussian random fields.  In this context, we consider the asymptotic of the extremal shot noise in the case when the weight distribution $G$ is regularly varying and the intensity of points goes to infinity.  

Several recent works also put the emphasis on statistical properties of max-stable random fields. The dependence structure
is investigated thanks to the extremal coefficient \cite{L83,SchTaw} or the extremal index \cite{S06}.

The structure of the paper is the following: in Section \ref{sec:2} we study general properties of extremal shot noises such as boundedness, regularity, ergodicity. In Section \ref{sec:3}, we consider extremal shot noise with regularly varying heavy tailed-weight distribution and prove convergence of the rescaled extremal shot noise to a max-stable random field when the intensity goes to infinity. Our approach is based on a functional point processes approach that is a functional version of the so-called Peak Over Threshold method in extreme value theory. Section \ref{sec:4} is devoted to the limit max-stable extremal shot noise: we give explicit formulas for its extremal coefficient, its extremal index  and we also consider the extremal points of the underlying Poisson point process.

\section{Properties of the extremal shot noise}\label{sec:2}
Let $(\Omega,\cF,\bbP)$ be a probability space. We consider $\Phi$ a Poisson point process on $\bbR^d\times(0,+\infty)$ with intensity $dx G(dm)$, where $G$ denotes a measure on $(0,+\infty)$. We suppose that the tail function $\bar G$ defined by $\bar G(u)=G((u,+\infty))$ is finite for any $u>0$. A generic point of the point process $\Phi$ is denoted by $\phi=(x_\phi,m_\phi)$. Denote by $\cF_\Phi$ the $\sigma$-field generated by $\Phi$, i.e. the $\sigma$-field generated by the random variables ${\rm card}(\Phi\cap A)$, for all $A$ compact set of $\bbR^d\times(0,+\infty)$. We suppose that $\Phi$ is $\cF$-measurable, i.e. $\cF_\Phi\subset \cF$.
The state space for the Poisson point process $\Phi$ is  the set $M_p$ of locally finite subsets of $\bbR^d\times (0,+\infty)$ endowed with the $\sigma$-field $\cM_p$ generated by the applications $M_p\to \bbN, N\mapsto {\rm card}(N\cap A)$,  $A$ compact set in $\bbR^d\times (0,+\infty)$. We denote by ${\rm PPP}(dx G(dm))$ the law of $\Phi$, where PPP stands for Poisson point process.

Let $h:\bbR^d\to[0,+\infty]$ be a measurable function and for $\phi\in\Phi$, denote by $h_\phi$ the function $h_\phi(y)=m_\phi h(y-x_\phi)$. We consider the extremal shot noise $M$ on $\bbR^d$ defined by 
\begin{equation}\label{eq:defshotnoise}
M(y)=\sup\{ h_\phi(y)\ ;\ \phi\in\Phi\},\quad y\in\bbR^d,
\end{equation}
where we stress that the supremum may be equal to $+\infty$. The state space for the extremal shot noise $M$ is the product space $[0,+\infty]^{\bbsR^d}$  endowed with the product $\sigma$-field generated by the projections $[0,+\infty]^{\bbsR^d}\to [0,+\infty],\ f\mapsto f(y)$, for $y\in\bbR^d$. We denote by ${\rm ESN}(h,G)$ the law of $M$, where ESN stands for extremal shot noise.

We first state a simple measurability property of the extremal shot-noise:
\begin{prop}\label{prop:meas} Let $\Phi\sim{\rm PPP}(dxG(dm)$ and $M\sim {\rm ESN}(h,G)$ the associated extremal shot noise. The application 
$$M:\left\{\begin{array}{ccc} (\Omega\times\bbR^d,\cF_\Phi\otimes \cB(\bbR^d))&\to& ([0,+\infty],\cB([0,+\infty]))\\ (\omega,y)&\mapsto& M(\omega,y)\end{array}\right.$$
is measurable. 
\end{prop}
{\it Proof:}
We introduce a measurable enumeration of the points of the Poisson point process $\Phi$ (cf \cite{DVJ}). There is a family of measurable applications $\phi_i: (\Omega,\cF_\Phi)\to (\bbR^d\times (0,+\infty),\cB(\bbR^d\times (0,+\infty)))$, $i\geq 1$, such that $\Phi=\{\phi_i,i\geq 1\}$. Then the extremal shot noise $M$ is given by 
$M=\sup_{i\geq 1} h_{\phi_i}$. For each $i\geq 1$, the application $(\omega,y)\mapsto h_{\phi_i(\omega)}(y)$ is $\cF_\Phi\otimes \cB(\bbR^d)-\cB([0,+\infty])$ measurable. The countable pointwise supremum is also measurable.\CQFD

\begin{remark}\label{rk:1} A closely related model of extremal shot noise is given by
$$\tilde M(y)=\sup\{m_{\tilde\phi}+ \tilde h(y-x_{\tilde\phi})\ ;\ \tilde\phi\in\tilde\Phi\}\in [-\infty,+\infty],$$
with  $\tilde h:\bbR^d\to [-\infty,+\infty]$ a measurable function, $\tilde G$ a measure on $(-\infty,+\infty)$ such that $\tilde G((u,+\infty))$ for all $u>0$ and $\tilde\Phi$ a Poisson point process on $\bbR^d\times\bbR$ with intensity $ dx\times \tilde G(dm)$. Introduce the change of variables $M(y)=\exp(\tilde M(y))$, $h=\exp(\tilde h)$ and $\Phi=T(\tilde\Phi)$ where $T:(x,m)\mapsto (x,e^m)$. An application of the mapping Theorem (see e.g. \cite{King}) shows that $\Phi\sim {\rm PPP}(dx G(dm))$  where $G$ has tail function  $\bar G(u)=\tilde G((\ln u,+\infty))$; hence $M\sim {\rm ESN}(h,G)$. 
\end{remark}

\subsection{Invariance properties}
A first basic feature of the extremal shot noise is stationarity:
\begin{prop}\label{prop:stat} Let $M\sim {\rm ESN}(h,G)$. Then $M$ is a stationary random field; i.e.
$$\forall z\in\bbR^d,\quad M(\cdot-z)\stackrel{\cL}{=} M(\cdot).$$ 
\end{prop}
{\it Proof:} 
Since the Lebesgue measure on $\bbR^d$ is invariant under translation, the Poisson Point Process $\Phi\sim {\rm PPP}(dxG(dm))$ is invariant under the transformation of $\bbR^d\times (0,+\infty)$ defined by $(x,m)\mapsto (x+z,m)$. The translated point process $\Phi+(z,0)$ has hence distribution ${\rm PPP}(dxG(dm))$. Now the translated extremal shot noise $M(\cdot+z)$ is based on the the extremal shot noise based on the translated point process $\Phi+(z,0)$ and hence has distribution $M\sim {\rm ESN}(h,G)$.\CQFD
\\ \ \\
Stationarity is not always a desirable feature in applications. Note that non-stationary models can be designed by replacing the Lebesgue measure $dx$ by a general $\sigma$-finite measure. Most of our results below would still hold true.

The class of extremal shot noises enjoys a nice property of stability  with respect to pointwise maximum; this is closely related to the property of max-infinitely divisibility (see \cite{EHV90}).
\begin{prop} 
\begin{enumerate}
\item Let $M_i,\ 1\leq i\leq n$, be independent extremal shot noises with distribution ${\rm ESN}(h,G_i)$
respectively. Then $M= \vee_{i=1}^n M_i$ has distribution ${\rm ESN}(h,G)$, with  $G=\sum_{i=1}^n G_i$.
\item Let $M\sim {\rm ESN}(h,G)$. Then $M$ is max-infinitely divisible: for all $n\geq 1$, $M\stackrel{\cL}{=}\vee_{i=1}^n M_i$ where $M_i,\ 1\leq i\leq n,$ are i.i.d. random fields with distribution ${\rm ESN}(h,n^{-1}G)$.
\end{enumerate}
\end{prop}
{\it Proof:} 
We prove only the first item, the second is a straightforward consequence. Let $\Phi_i, 1\leq i\leq n,$ be independent Poisson processes with distribution ${\rm PPP}(dxG_i(dm))$ respectively. The associated shot noise are 
$M_i=\sup_{\phi\in\Phi_i} h_\phi$, $1\leq i\leq n$. Denote by $\Phi=\bigcup_{i=1}^n \Phi_i$. From the superposition Theorem (see e.g. \cite{King}), $\Phi$ is a Poisson point process with intensity $G=\sum_{i=1}^n G_i$. Furthermore, 
$M=\vee_{i=1}^n M_i$ is the extremal shot noise associated with the Poisson Point Process $\Phi$, i.e. $M=\sup_{\phi\in\Phi} h_\phi$. This implies $M\sim {\rm ESN}(h,G)$.\CQFD

\subsection{Finite-dimensional distributions}
We give a condition for the extremal shot noise to be finite and characterize its finite dimensional distributions. 
An important quantity is the coefficient $\alpha(h,G)$ defined by
$$\alpha(h,G)=\inf \left\{u>0; \int_{\bbsR^d}\bar G\left(u/h(x)\right)dx<+\infty\right\},$$
with the convention $\alpha(h,G)=+\infty$ if the set is empty,  $\bar G\left(u/h(x)\right)=0$ if $h(x)=0$ and $\bar G\left(u/h(x)\right)=\bar G(0)=G((0,+\infty])$ if $h(x)=+\infty$. As we will see below, the coefficient $\alpha(h,G)$ appears to be the left end-point of the support of the distribution of $M(y)$.\\

\begin{prop}\label{prop:pr1} Let $\Phi\sim {\rm PPP}(dxG(dm))$ and $M\sim {\rm ESN}(h,G)$ the associated extremal shot noise.
\begin{enumerate} 
\item The cumulative distribution function (cdf) of $M(y)$ is given by
$$\bbP(M(y)\leq u)=\left\{\begin{array}{ll}
0 & \mbox{if} \  u<\alpha(h,G), \\
\exp\left(-\int_{\bbsR^d}\bar G\left(u/h(x)\right)dx\right)& \mbox{otherwise}.\\
\end{array}\right.$$
More generally, the multivariate cdf is given by: for $y_1,\ldots,y_k\in \bbR^d$ and $u_1,\ldots,u_k\in\bbR$,
\begin{eqnarray*}
&&\bbP(M(y_1)\leq u_1,\ldots,M(y_k)\leq u_k)\\
&=&\left\{\begin{array}{ll}
0 & \mbox{if} \  \wedge_{1\leq i\leq k}u_i<\alpha(h,G), \\
\exp\left(-\int_{\bbsR^d}\bar G\left(\wedge_{1\leq i\leq k}\{u_i/h(y_i-x)\}\right)dx\right)& \mbox{otherwise}.\\
\end{array}\right.
\end{eqnarray*}
\item If $h$ is finite almost everywhere, the following zero/one law holds: for all $y\in\bbR^d$,
$$\bbP(M(y)=+\infty)=\left\{\begin{array}{l} 0 \ \mbox{if}\ \alpha(h,G)<+\infty \\ 1 \ \mbox{if}\ \alpha(h,G)=+\infty\end{array}\right.,\quad y\in\bbR^d.$$  
\end{enumerate}
\end{prop}
{\it Proof of Proposition \ref{prop:pr1}:}
The event $\{M(y)\leq u\}$ can be written as 
\begin{eqnarray*}
\{M(y)\leq u\}&=&\left\{\forall \phi\in\Phi\ ,\ m_\phi h(y-x_\phi)\leq u \right\}\\
&=&\left\{\Phi\cap A=\emptyset\right\}
\end{eqnarray*}
with $A=\{(x,m)\in\bbR^d\times(0,+\infty); mh(y-x)>u\}$. This shows that $M(y)$ is $\cF_\Phi$-measurable. Using the avoidance probability for the Poisson Point Process $\Phi$,
\begin{eqnarray}
\bbP(M(y)\leq u)&=&\exp\left(-\int_{\bbsR^d\times(0,+\infty)}1_{\{m>u/h(y-x)\}} dxG(dm)\right)\nonumber\\
&=&\exp\left(-\int_{\bbsR^d}\bar G\left(u/h(x)\right)dx\right) \label{eq:cdf}.
\end{eqnarray}
This gives the cdf of $M(y)$. The multivariate cdf is obtained in a similar way:
\begin{eqnarray}
&&\bbP(M(y_1)\leq u_1,\ldots,M(y_k)\leq u_k)\nonumber\\
&=&\bbP\left(\forall \phi\in\Phi,\ \forall 1\leq i\leq k,\ m_\phi h(y_i-x_\phi)\leq u_i \right)\nonumber\\
&=&\bbP\left(\Phi\cap \left\{(x,m)\in\bbR^d\times\bbR^+; m>\wedge_{1\leq i\leq k}\{u_i/h(y_i-x)\}\right\}=\emptyset\right)\nonumber\\
&=&\exp\left(-\int_{\bbsR^d\times(0,+\infty)}1_{\left\{m>\wedge_{1\leq i\leq k}\{u_i/h(y_i-x)\}\right\}} dxG(dm)\right)\nonumber\\
&=&\exp\left(-\int_{\bbsR^d}\bar G\left(\wedge_{1\leq i\leq k}\{u_i/h(y_i-x)\}\right)dx\right)\nonumber.
\end{eqnarray}
We now prove the zero/one law. If $\alpha(h,G)=+\infty$,  $\int_{\bbsR^d}\bar G\left(u/h(x)\right)dx=+\infty$ for all $u>0$. This implies $\bbP(M(y)\leq u)=0$ for all $u>0$ and $\bbP(M(y)=+\infty)=1$.
In the case $\alpha(h,G)<+\infty$,  $\int_{\bbsR^d}\bar G\left(u/h(x)\right)dx<+\infty$ for all $u>\alpha(h,G)$  and 
$$\bbP(M(y)\leq u)=\exp\left(-\int_{\bbsR^d}\bar G\left(u/h(x)\right)dx\right)>0.$$
For all $x\in\bbR^d$ such that $h(x)<+\infty$, the nonincreasing convergence $\bar G\left(u/h(x)\right)\to 0$ holds as $u\to+\infty$. If $h$ is finite almost everywhere, then the monotone convergence Theorem we obtain
$$\bbP(M(y)<+\infty)=\lim_{u\to+\infty}\bbP(M(y)<u)=1.$$
This proves the zero/one law.\CQFD

\begin{example}\label{ex:1} Consider the case when for all $u>0$, $\bar G(u)=u^{-\xi}$ for some $\xi>0$ and $\int_{\bbsR^d} h(x)^\xi dx <+\infty$. Then $\alpha(h,G)=0$ and we recover the stationary max-stable model from Smith \cite{Sm} and Schlather \cite{Schla}: the univariate distribution of the extremal shot noise is a Fréchet distribution with cdf
$$\bbP(M(y)\leq u)=\exp\left(-u^{-\xi}\int_{\bbsR^d}h(x)^\xi dx \right);$$
and the multivariate distribution is given by
$$\bbP(M(y_1)\leq u_1,\ldots,M(y_k)\leq u_k)=\exp\left(-\int_{\bbsR^d}\left(\wedge_{1\leq i\leq k}\{h(y_i-x)^\xi u_i^{-\xi}\}\right)dx\right).$$
Moreover, it satisfies the max-stability functional equation: for all $\theta>0$ 
$$\bbP(M(y_1)\leq u_1,\ldots,M(y_k)\leq u_k )^\theta= \bbP(M(y_1)\leq \theta^{-1/\xi}u_1,\ldots,M(y_k)\leq \theta^{-1/\xi}u_k).$$
See Section \ref{sec:4} for further discussions on the max-stable case.
\end{example}

\begin{example}\label{ex:2} Consider the case when $G$ has an exponential distribution, i.e. $\bar G(u)=e^{-u}$ for all $u>0$. Suppose that the function $h$  is bounded and satisfies for large enough $x$:
$$\frac{\gamma^-(|x|)}{\ln(|x|)}\leq h(x)\leq \frac{\gamma^+(|x|)}{\ln(|x|)}$$
with 
$$\lim_{r\to +\infty}\gamma^-(r)=\lim_{r\to +\infty}\gamma^+(r)=\gamma\in [0,+\infty].$$
For $u>0$, the function $x\mapsto \bar G\left(u/h(x)\right)$ is locally bounded and satisfies for large $x$,
$$|x|^{-\frac{u}{\gamma^-(|x|)}} \leq \bar G\left(u/h(x)\right)\leq |x|^{-\frac{u}{\gamma^+(|x|)}}.$$
Since $x\mapsto |x|^{-\frac{u}{\gamma}}$ is integrable near infinity if and only if $\frac{u}{\gamma}>d$, we obtain $\alpha(h,G)=\gamma d$.
\end{example}

\begin{remark} Another interpretation of the coefficient $\alpha(h,G)$ is the contribution from points at infinity: it can be shown that if $h$ is locally bounded or $G$ is finite,  then for all $y\in\bbR^d$:
$$\sup\{m_\phi h(y-x_\phi)\ ;\ \phi\in\Phi \ ,\ |x_\phi|\geq R\}\to \alpha(h,G)$$
almost surely as $R\to+\infty$.
\end{remark}

%The function $R\mapsto M_R(y)$ is nonincreasing and nonnegative, and hence converges almost surely  as $R\to+\infty$.
%A similar computation as in equation (\ref{eq:cdf}) yields
%$$\bbP(M^R(y)\leq u)=\exp\left(-\int_{\bbsR^d}\bar G(u/h(y-x))1_{\{|x|\geq R\}}dx\right).$$
%We easily see that if $h$ is locally bounded or $G$ is finite,
%$$\int_{\bbsR^d}\bar G(u/h(y-x))1_{\{|x|\geq R\}}dx <+\infty \quad \mbox{\ if\ and\ only\ if\ }\quad \int_{\bbsR^d}\bar G(u/h(x))dx <+\infty.$$ 
%As a consequence, if $u<\alpha$, $\bbP(M^R(y)<u)=0$  and if $u>\alpha$, monotone convergence implies
%$$\lim_{R\to+\infty}\bbP(M^R(y)\leq u)=\lim_{R\to+\infty} \exp\left(-\int_{\bbsR^d}\bar G(u/h(y-x))1_{\{|x|\geq R\}}dx\right)=1.$$
%This proves the almost sure convergence $M^R(y)\to\alpha$ as $R\to+\infty$.

\subsection{Boundedness}
We now explore the path properties of the random field $M$ and consider first the boundedness property.
Let $M\sim{\rm ESN}(h,G)$ and $A$ be a nonempty measurable subset of $\bbR^d$. Define
$$M_{A}(y)=\sup_{z\in A} M(y+z)\in [0,+\infty].$$ 
We easily see that $M_{A}\sim {\rm ESN}(h_A,G)$ where 
$$h_A(x)=\sup_{z\in A} h(x+z),\quad x\in\bbR^d.$$
An application of Proposition \ref{prop:pr1} yields the following interesting corollary:
\begin{corol}\label{cor1} Let $M\sim {\rm ESN}(h,G)$ and $A$ be a non-empty bounded measurable subset of $\bbR^d$.
\begin{enumerate} 
%\item The random variable $\sup_{z\in A} M(z)$ is $\cF_{\Phi}$-measurable.
\item The distribution of $\sup_{z\in A} M(z)$ has cdf
$$\bbP\left(\sup_{z\in A} M(z)\leq u\right)=\left\{\begin{array}{ll}
0 & \mbox{if} \  u<\alpha(h_A,G), \\
\exp\left(-\int_{\bbsR^d}\bar G\left(u/h_A(x)\right)dx\right)& \mbox{if }\  u\geq \alpha(h_A,G).\\
\end{array}\right.$$
\item If $A$ is bounded, $h$ is locally bounded and $\alpha(h_A,G)<+\infty$, then $\sup_{z\in A} M(z)<+\infty$ almost surely.
\end{enumerate}
\end{corol}
The coefficient $\alpha(h_A,G)$ is the left end-point of the support of the distribution of $\sup_{z\in A} M(z)$.

\begin{example}\label{ex:3}
The inequality $\alpha(h,G)\leq \alpha(h_A,G)$ always holds true. But it can be strict as the following example shows: let $d=1$, $A=[-1,1]$, let $G$ be the exponential distribution i.e. $\bar G(u)=e^{-u}$ for $u>0$, and let  
$$h(x)=\sum_{n\geq 1}\frac{\gamma}{1+\ln(1+|x|)}1_{\{n< |x| < n+n^{-2} \}}$$
for some $\gamma>0$. Since $\sum_{n\geq 1}n^{-2}<+\infty$, we easily see that $\int_{\bbsR}1_{\{h(x)>0\}}dx<+\infty$ and this implies $\alpha(h,G)=0$.
On the other hand, the function $h_A$ satisfies for all $|x|>1$
$$\frac{\gamma}{1+\ln(2+|x|)} \leq h_A(x)\leq \frac{\gamma}{1+\ln(|x|)}.$$
From  example \ref{ex:2} above, we conclude that $\alpha(h_A,G)=\gamma$.
\end{example}

\begin{remark}\label{rk:3} 
The following relations are worth noting:\\
- if $A_1\subset A_2$, then $\alpha(h_{A_1},G)\leq \alpha(h_{A_2},G)$;\\
- if $A_2=A_1+x$, then $\alpha(h_{A_1},G)=\alpha(h_{A_2},G)$;\\
- if $A_3=A_1\cup A_2$, then $\alpha(h_{A_3},G)=\alpha(h_{A_1},G)\vee \alpha(h_{A_2},G)$.\\
The last statement is proved as follows:  $h_{A_3}=h_{A_1}\vee h_{A_2}$ so that for all $u>0$,
$$\bar G(u/h_{A_i}(x)) \leq \bar G(u/h_{A_3}(x))\leq  \bar G(u/h_{A_1}(x))+\bar G(u/h_{A_2}(x)),\quad x\in\bbR^d,i\in\{1,2\};$$
this implies that the integral $\int\bar G(u/h_{A_3}(x))dx$ is finite if and only if $\int\bar G(u/h_{A_i}(x))dx$ is finite for $i=1$ and $i=2$; as a consequence, $\alpha(h_{A_3},G)=\alpha(h_{A_1},G)\vee \alpha(h_{A_2},G)$.\\
Using these three properties, one can show that for any  bounded sets $A_1$, $A_2$ with nonempty interior, we have $\alpha(h_{A_1},G)=\alpha(h_{A_2},G)$. The common value is denoted by $\alpha^+(h,G)$.
\end{remark}

Of particular interest is the property that the extremal shot noise is bounded on compact sets. Define the function $h^+$ by
\begin{equation}\label{eq:defh+}
h^+(x)=\sup\{h(x+z),|z|\leq 1\},\quad x\in\bbR^d,
\end{equation}
and let $\alpha^+(h,G)=\alpha(h^+,G)$ (this is consistant with the previous definition of $\alpha^+(h,G)$).
\begin{corol}\label{cor2}
Suppose $h$ is locally bounded and  $M\sim{\rm ESN}(h,G)$. 
Then $M$ is almost surely bounded on compact sets if and only if $\alpha^+(h,G)<+\infty$.
\end{corol}
{\it Proof:} Denote by $(K_n)_{n\geq 1}$ an increasing sequence of compact sets with non-empty interiors and such that $\bbR^d=\cup_{n\geq 1}{\rm int}(K_n)$. As stated in Remark \ref{rk:3}, for all $n\geq 1$ $\alpha(h_{K_n},G)=\alpha^+(h,G)$. From Corollary \ref{cor1} and using the zero/one law, $M$ is a.s. bounded on $K_n$ if and only if $\alpha(h_{K_n},G)<+\infty$. This implies that $M$ is a.s. bounded on all the compact sets $(K_n)_{n\geq 1}$ if and only if $\alpha^+(h,G)<+\infty$. Finally, since for all compact $K$ there exists $n\geq 1$ such that $K\subseteq K_n$, we see that $M$ is a.s. bounded on compact sets if and only if $\alpha^+(h,G)<+\infty$. \CQFD

\begin{example} 
We provide an example where $\alpha(h,G)=0$ and $\alpha^+(h,G)=+\infty$. This implies that the corresponding extremal shot noise $M$ is almost surely finite at all point, but unbounded from above on any open set. Let $d=1$, $A=[-1,1]$, $\bar G(u)=e^{-u}$ for $u>0$ and
$$h(x)=\sum_{n\geq 1}\frac{\gamma(|x|)}{1+\ln(1+|x|)}1_{\{n< |x| < n+n^{-2} \}},\quad x\in\bbR,$$
with $\gamma:[0,+\infty)\to [0,+\infty)$ a nondecreasing function such that $\gamma(u)\to+\infty$ and $\gamma(u)/\ln(u)\to 0$ as $u\to+\infty$. Since $h$ is bounded and the set $\{h>0\}$ has finite Lebesgue measure, $\alpha(h,G)=0$. On the other hand, if $|x|>1$,
$$ h_A(x)\geq \frac{\gamma(|x|-1)}{1+\ln(2+|x|)}.$$
From Example \ref{ex:2} above, we conclude that $\alpha^+(h,G)=+\infty$.
\end{example}

\subsection{Regularity}
We proceed with regularity properties of the extremal shot noises. The regularity of the extremal shot noise $M\sim{\rm ESN}(h,G)$ depends on the regularity of the shape function $h$.
We need the following definition: for all $\varepsilon>0$, define
$$h^-_{\varepsilon}(x)=\inf\{h(x+z); |z|\leq\varepsilon\},\quad x\in\bbR^d,$$
and $\alpha^-(h,G)=\lim_{\varepsilon\to 0} \alpha(h^-_{\varepsilon},G)$. 
Note that the following inequality always holds true: $0\leq \alpha^-(h,G)\leq \alpha(h,G)\leq \alpha^+(h,G)$.
\begin{prop}\label{prop:pr2} Let $M\sim{\rm ESN}(h,G)$.
\begin{enumerate}
\item If $h$ is lower semi-continuous, then $M$  is lower semi-continuous.
\item If $\alpha^-(h,G)=\alpha^+(h,G)$ and $h$ is upper semi-continuous, then $M$ is almost surely upper semi-continuous.
\item If $\alpha^-(h,G)=\alpha^+(h,G)$ and $h$ is continuous, then  $M$ is almost surely continuous.
\end{enumerate}
\end{prop}

\begin{remark} It is worth noting that a necessary condition for the random field $M$ to be upper semi-continuous is that 
$\alpha(h,G)=\alpha^+(h,G)$. Otherwise, if $\alpha(h,G)<\alpha^+(h,G)$, Proposition \ref{prop:pr1} and Corollary \ref{cor1},  imply that with positive probability $\alpha(h,G)\leq M(0)<\alpha^+(h,G)$  whereas for all $\varepsilon>0$, $\sup\{M(y); |y|\leq \varepsilon\}\geq\alpha^+(h,G)$  with probability $1$; this implies that the random field $M$ cannot be upper semi-continuous with probability $1$. 
In Proposition \ref{prop:pr2}, we prove upper semi-continuity under the stronger condition $\alpha^{-}(h,G)=\alpha^+(h,G)$;  this condition might be stronger than necessary but is required in our proof to give a lower bound for the random field $M$ uniform on compact sets.
\end{remark}

The proof of Proposition \ref{prop:pr2} relies on the following lemma that gives some insight into the structure of the extremal shot noise:
\begin{lemme}\label{lem:lem2}
Let $M\sim{\rm ESN}(h,G)$. 
\begin{enumerate}
\item For all $y\in\bbR^d$, $M(y)\geq \alpha^{-}(h,G)$;
\item For all compact $K\subset\bbR^d$ and $u>\alpha^+(h,G)$, there exists a (random) finite subset $\Phi_{K,u}\subseteq \Phi$ such that
\begin{equation}\label{eq:repres}
M(y)\vee u=\max\{h_\phi(y)\ ;\ \phi\in\Phi_{K,u}\}\vee u,\quad \mbox{ for\ all\ }y\in K.
\end{equation}
\end{enumerate}
\end{lemme}

We first show how Lemma \ref{lem:lem2} implies Proposition \ref{prop:pr2} and then proceed with the proof of Lemma \ref{lem:lem2}.\\
{\it Proof of Proposition \ref{prop:pr2}:}
A supremum of lower semi-continuous functions is always lower semi-continuous; this implies the first point.
The third point is a consequence of the two first points since a function is continuous if and only if it is both upper and lower semi-continuous. It remains to prove the second point. Suppose $h$ is upper semi-continuous and $\alpha^{-}(h,G)=\alpha^+(h,G)$. Let $K\subseteq \bbR^d$ be a compact set. Lemma \ref{lem:lem2} implies that for all $u>\alpha^+(h,G)$, the random field $M\vee u$ is almost surely upper semi-continuous on $K$ as a maximum of a finite number of upper semi-continuous functions. Let $u_n\to \alpha^+(h,G)$ be a decreasing sequence. The sequence of upper semi-continuous random fields $M\vee u_n$ converges uniformly on $K$ to $M\vee \alpha^+(h,G)$ and hence $M\vee \alpha^+(h,G)$ is almost surely upper semi-continuous on $K$. The compact $K$ being arbitrary, $M\vee \alpha^+(h,G)$ is almost surely upper semi-continuous on $\bbR^d$. Since from Lemma \ref{lem:lem2}, $M\geq  \alpha^-(h,G)$, the conclusion comes from the condition  $\alpha^-(h,G)=\alpha^+(h,G)$. \CQFD\\
{\it Proof of Lemma \ref{lem:lem2}:} 
Let $\varepsilon>0$ and $(y_n)_{n\geq 1}$ be a sequence of points in $\bbR^d$ such that $\bigcup_{n\geq 1} \{y_i+z;\ |z|\leq \varepsilon\}=\bbR^d$. Note that, for all $n\geq 1$,
$$\inf\{M(y_n+z); |z|\leq \varepsilon\}\geq \sup\{m_\phi h^{-}_{\varepsilon}(y_n-x_\phi); \phi\in\Phi \}.$$
Proposition \ref{prop:pr1} gives the distribution of the extremal shot noise ${\rm ESN}(h^-_{\varepsilon},G)$ and implies that for all $n\geq 1$, $\sup\{m_\phi h^{-}_{\varepsilon}(y_n-x_\phi); \phi\in\Phi \}\geq \alpha(h^-_{\varepsilon},G)$ almost surely. Hence, almost surely $M(y)\geq \alpha(h^-_{\varepsilon},G)$ for all $y\in \bbR^d$. The result is proved letting $\varepsilon\to 0$.\\
To prove the second point, define $\Phi_{K,u}=\{\phi\in\Phi;\ \sup_{y\in K} h_\phi(y)>u\}$. This definition ensures that equation \eqref{eq:repres} is satisfied. It remains to verify that $\Phi_{K,u}$ is finite almost surely. This is the case because the random variable ${\rm card}(\Phi_{K,u})$ has a Poisson distribution with mean
$$\int_{\bbsR^d\times (0,+\infty)} 1_{\{\sup_{y\in K} mh(y-x)>u\}} G(dm)dx=\int_{\bbsR^d}\bar G(u/h_K(x))dx$$ 
which is finite since $u>\alpha^{+}(h,G)\geq \alpha(h_K,G)$.\CQFD

According to Proposition \ref{prop:pr2}, the condition $\alpha^-(h,G)=\alpha^+(h,G)$ plays an important role. We give now necessary conditions so that it holds true.  We disregard the case when $\alpha^-(h,G)=\alpha^+(h,G)=+\infty$ because in this case the random field $M$ is almost surely constant and equal to $+\infty$. Note that condition $\alpha^+(h,G)<+\infty$ implies that $G$ is finite or $h$ is bounded. 
\begin{prop}\label{prop:alpha}
Suppose $G$ is finite or $h$ is bounded.  
\begin{enumerate}
\item If there exist $\gamma>d$ and $C>0$ such that  
\begin{equation}\label{eq:alpha1}
\lim_{x\to\infty} h(x)\bar G^{\leftarrow}( C \|x\|^{-\gamma})=0 \quad {\rm with}\quad \bar G^{\leftarrow}(t)=\inf\{u>0; \bar G(u)\leq t\},
\end{equation}
then $\alpha^-(h,G)=\alpha^+(h,G)=0$.
\item Let $\varepsilon>0$. If for all $\delta\in(0,1)$, there exists  $R>0$ such that for all $\|x\|>R$ and $\|z\|\leq \varepsilon$ 
\begin{equation}\label{eq:alpha2}
(1-\delta)h(x)\leq h(x+z) \leq (1+\delta)h(x),
\end{equation}
then   $\alpha^-(h,G)=\alpha^+(h,G)<+\infty$.
\end{enumerate}
\end{prop}
{\it Proof of Proposition \ref{prop:alpha}:}
For the first point, it is enough to check that 
$${\rm for\ all\ } u>0,\quad \int_{\bbsR^d} \bar G(u/h^+(x))dx <+\infty.$$
The condition $G$ finite or $h$ bounded ensures that the function  $x\mapsto \bar G(u/h^+(x))$ is locally integrable.
The function $\bar G^{\leftarrow}$ satisfies $u\geq \bar G^{\leftarrow}(t)$ if and only if $\bar G(u)\leq t$. Equation \eqref{eq:alpha1} ensures that for large $x$, $h(x)\bar G^{\leftarrow}( C \|x\|^{-\gamma})\leq u$ so that
$h^+(x)\bar G^{\leftarrow}( C (\|x\|-1)^{-\gamma})\leq u$ and 
$\bar G(u/h^+(x))\leq C(\|x\|-1)^{-\gamma}$. As a consequence, the function  $x\mapsto \bar G(u/h^+(x))$ is integrable at infinity. This proves the first point.\\
For the second point, notice that equation \eqref{eq:alpha2} implies that the functions $h^-_\varepsilon$ and $h^+_\varepsilon$ are equivalent at infinity, where
$$h^+_\varepsilon(x)=\sup\{h(x+z); \|z\|\leq\varepsilon\}\quad {\rm and}\quad  h^-_\varepsilon(x)=\inf\{h(x+z); \|z\|\leq\varepsilon\}.$$
This implies that $\alpha(h^-_\varepsilon,G)=\alpha(h^+_\varepsilon,G)$ and hence $\alpha(h^-,G)=\alpha(h^+,G)$.\CQFD

\subsection{Separability}
We consider the separability property of extremal shot noises (see \cite{Bill}). Let $D$ be a countable dense subset of $\bbR^d$. A function $f:\bbR^d\to [0,+\infty]$ is said to be $D$-separable if for all $t\in\bbR^d$, there is a sequence $(t_n)_{n\geq 1}$ of points of $D$ such that $t_n\to t$ and $f(t_n)\to f(t)$. We say that $f$ is universally separable (US) if $f$ is $D$-separable for all dense countable subset $D$ of $\bbR^d$. For example, any continuous function is US since $f(x)=\lim f(x')$ where the limit is taken when $x'\to x,x'\in D$. 
In the framework of extremal shot noises, the following observation plays an important role: the class of universally separable upper semi-continuous (USUSC) functions is closed for the topology of uniform convergence on compact sets, is translation invariant, and is stable under finite pointwise maximum. 

\begin{prop}\label{prop:pr2bis} Let $M\sim{\rm ESN}(h,G)$. Suppose $h$ is USUSC and $\alpha^{-}(h,G)=\alpha^{+}(h,G)$. Then $M$ is almost surely USUSC.
\end{prop}
{\it Proof:}
 Let $O$ be an open and relatively compact subset of $\bbR^d$, and denote by $K$ its closure. According to Lemma \ref{lem:lem2}, for all $u>\alpha^{+}(h,G)$, there exists a finite set $\Phi_{K,u}$ such that
$$M\vee u=\sup\{h_\phi;\ h\in\Phi_{K,u}\}\vee u\quad {\rm on\ }K.$$
Since $\Phi_{K,u}$ is finite and the functions $h_\phi$ are  USUSC, $M\vee u$ is USUSC. With a similar argument as in the proof of Proposition \ref{prop:pr2}, the assumption $\alpha^{-}(h,G)=\alpha^{+}(h,G)$ implies that $M\vee u$ converge uniformly to $M$ as $u\to \alpha^{+}(h,G)$. Hence $M$ is also USUSC. \CQFD

\begin{remark}
The class of US functions is not stable under pointwise maximum as the following example shows: $h_1(x)=1_{\{x>0\}}$ and $h_2(x)=1_{\{x<0\}}$ are US (and also lower semi-continuous), but $h_1\vee h_2(x)=1_{\{x\neq 0\}}$ is not US. This is the reason why we need to consider the class USUSC in Proposition \ref{prop:pr2bis}.
\end{remark}

\subsection{Ergodicity}
We end this section with the mixing properties of the extremal shot noise (see \cite{Bill}).  Note that related results for max-stable process are obtained by Stoev \cite{S08}, using the spectral representation in terms of extremal integrals. 
The stationary random field $M$ is said to be $\alpha$-mixing in $[0,+\infty]^{\bbsR^d}$
if for any $A,B\in\cB([0,+\infty]^{\bbsR^d})$,
\begin{equation}\label{eq:mix}
\lim_{v\to\infty} \bbP(M\in A \ ,\ \tau_v M\in B)= \bbP(M\in A)\bbP(M\in B),
\end{equation}
where $\tau_vM(\cdot)=M(\cdot+v)$, $v\in\bbR^d$.
When the random field $M$ has almost surely continuous paths, it can be considered as a random element of $\bbC(\bbR^d,[0,+\infty])$ endowed with the metric of uniform convergence on compact sets and the Borelian $\sigma$-field $\cB(\bbC(\bbR^d,[0,+\infty]))$. In this case, we say that the random field $M$ is $\alpha$-mixing in $\bbC(\bbR^d,[0,+\infty])$ if \eqref{eq:mix} holds true for all $A,B\in\cB(\bbC(\bbR^d,[0,+\infty]))$.
\begin{prop}\label{prop:pr3} Let $M\sim {\rm ESN}(h,G)$.
\begin{enumerate} 
\item The extremal shot noise $M$ is $\alpha$-mixing in $[0,+\infty]^{\bbsR^d}$;
\item If $h$ is continuous and $\alpha^{-}(h,G)=\alpha^+(h,G)$, then $M$ is $\alpha$-mixing in $\bbC(\bbR^d,[0,+\infty])$.
\end{enumerate}
\end{prop}
Recall that $\alpha$-mixing implies ergodicity. The mixing property  will be proved using Lemma \ref{lem:lem2} and the independence property of Poisson point process.

\noindent
{\it Proof of Proposition \ref{prop:pr3}:}
According to \cite{Bill}, it is enough to check the mixing property \eqref{eq:mix} for $A,B$ in a $\pi$-system generating $\cB([0,+\infty]^{\bbsR^d})$. Such a $\pi$-system is given by the finite intersections of sets of the form $\{f\in [0,+\infty]^{\bbsR^d}; f(y)\geq u\}$ for some $t\in\bbR^d$ and $u\geq 0$. Let $A,B$ be given by
\begin{eqnarray*}
A&=&\{f\in[0,+\infty]^{\bbsR^d}; f(y_1)\geq u_1,\ldots,f(y_k)\geq u_k\},\\
B&=&\{f\in[0,+\infty]^{\bbsR^d}; f(y_1')\geq u_1',\ldots,f(y_l')\geq u_l'\}.
\end{eqnarray*}
Since $M(t)\geq \alpha(h,G)$ almost surely, we can suppose w.l.o.g. that $u_i>\alpha(h,G), 1\leq i\leq k$ and 
$u_j'>\alpha(h,G), 1\leq j\leq l$. Let $K=\{y_i, y'_j; 1\leq i\leq k,\ 1\leq j\leq l\}$. 
From Lemma \ref{lem:lem2}, there is a finite point process $\Phi_{K,u}\subset\Phi$ such that equation (\ref{eq:repres}) holds for all $u>\alpha^+(h,G)$. Since $K$ is finite, a straightforward modification of Lemma \ref{lem:lem2} shows that its conclusion remains true in this case for $u>\alpha(h,G)$. 
In particular, let $u>\alpha(h,G)$ be the minimum of all $u_i$ and $u_j'$. The construction of $\Phi_{K,u}$ ensures that $\{M\in A\}=\{M\vee u\in A\}=\left\{\sup_{\phi\in\Phi_{K,u}}h_\phi  \in A\right\}$ almost surely, and the same result holds with $B$ replacing $A$. Since $\Phi_{K,u}$ is finite, for any $\varepsilon>0$, there is a compact $L\subset\bbR^d$ such that $\bbP(\Phi_{K,u} \subset L\cap (0,+\infty))\geq 1-\varepsilon$. Define $M_L=\sup\{h_\phi;\ \phi\in \Phi,\ x_\phi\in L\}$. Then $\bbP(M\vee u=M_L\vee u \mbox{\ on\ } K)\geq 1-\varepsilon$.
In the same way, we also have for any $v\in\bbR^d$, $\bbP(\tau_v M\vee u=\tau_v M_{L+v}\vee u \mbox{\ on\ } K)\geq 1-\varepsilon$, with $M_{L+v}=\sup\{h_\phi(y);\ \phi\in \Phi,\ x_\phi-v\in L\}$.
Hence we have,
$$ |\bbP(M_L\in A)- \bbP(M \in A)|\leq \varepsilon\ ,\ |\bbP(\tau_v M_{L+v}\in B)- \bbP(\tau_v M\in B)|\leq \varepsilon$$
and
$$|\bbP(M_L \in A\ ,\ \tau_vM_{L+v}\in B)- \bbP(M\in A\ ,\ \tau_vM\in B)|\leq 2\varepsilon.$$
The independence property of the Poisson Point Process $\Phi$ implies that $M_{L+v}$ and $M_L$ are independent for large $v$ since $L+v$ and $L$ are disjoint. Hence, for $v$ large enough,
$$\bbP(M_L \in A\ ,\ \tau_vM_{L+v}\in B)=\bbP(M_L \in A)\bbP(\tau_vM_{L+v}\in B).$$
Equation (\ref{eq:mix}) follows and this proves that  the extremal shot noise is $\alpha$-mixing in $[0,+\infty]^{\bbsR^d}$. 

According to Proposition \ref{prop:pr2}, the further conditions $h$ continuous and $\alpha^{-}(h,G)=\alpha^+(h,G)$ ensure that $M$ is almost surely continuous. The proof goes exactly the same way since the $\pi$-system consisting of the sets 
$$\{f\in\bbC(\bbR^d,[0,+\infty]);\ f(y_1)\geq u_1,\ldots,f(y_k)\geq u_k\}$$
generates the $\sigma$-field $\cB(\bbC(\bbR^d,[0,+\infty]))$.\CQFD\\

\section{Heavy-tailed extremal shot noises and their asymptotics}\label{sec:3}
In this section, we consider different asymptotics related to extremal shot noises when the weight measure $G$ is a probability measure with a regularly varying tail. We recall some facts about heavy-tailed probability measures, univariate extreme value theory  that will be useful in the sequel. For general references on this subject, see e.g. \cite{Res} or \cite{BGT}. 

We suppose that $G$ is a probability on $(0,+\infty)$ with tail function $\bar G\in RV_{-\xi},$ the set of functions regularly varying at infinity with exponent $-\xi<0$. This implies that $G$ belongs to the max-domain of attraction of the Fréchet distribution $F_\xi$ 
with cdf
$$F_\xi(x)= \exp(-x^{-\xi})1_{\{x>0\}}.$$
Indeed, there is a scaling $a_\lambda>0$ such that the distribution function $G$ satisfies 
\begin{equation}\label{eq:CVFrechet}
\lim_{\lambda\to+\infty}G(a_\lambda x)^\lambda=F_\xi(x).
\end{equation}
This has the following interpretation in terms of random variables: if $(X_i)_{i\geq 1}$ are i.i.d. with distribution $G$, then the renormalized maximum $a_n^{-1}\max_{1\leq i\leq n} X_i$ converges to the Fréchet distribution as $n\to+\infty$.
A possible choice for the renormalization function is 
\begin{equation}
a_\lambda=G^{\leftarrow}(1-\lambda^{-1}),\quad \lambda>0
\end{equation}
where $G^{\leftarrow}$ is the quantile function
$$G^{\leftarrow}(u)=\inf\{x>0; G((0,x])\leq u \},\quad 0<u<1.$$

Notice also that equation \eqref{eq:CVFrechet} implies the following estimate:
\begin{equation}\label{eq:asytail}
\lim_{\lambda\to+\infty}\lambda \bar G(a_\lambda x)=\bar G_\xi(x), \quad x>0,
\end{equation}
where $ \bar G_\xi(x)=x^{-\xi}$ and $G_\xi(dx)=\xi x^{-\xi-1}1_{\{x>0\}}dx$ is the corresponding infinite measure on $(0,+\infty)$. 

\subsection{Heavy-tailed extremal shot noise}
We consider the extremal shot noises $M_\lambda\sim {\rm ESN}(h,\lambda G)$ for some continuous shape function $h:\bbR^d\to [0,+\infty]$, intensity $\lambda>0$ and heavy-tailed probability $G$.

For $\xi>0$, we say that $h$ satisfies condition $({\bf C}_\xi)$ if:
\begin{enumerate}
%\item[$({\bf C_1})$\quad]$0<\liminf x^{\xi}\bar G(x)\leq  \limsup x^{\xi}\bar G(x) <+\infty$\ \ and\ \ $\int_{\bbsR^d} h^+(x)^\xi dx<+\infty$,
\item[$({\bf C}_\xi)$\quad] $\int_{\bbsR^d} h^+(x)^{\xi-\delta}dx<+\infty$ for some $\delta\in (0,\xi)$.
\end{enumerate}
where $h^+$ is given by \eqref{eq:defh+}. Condition $({\bf C}_\xi)$ implies that $h$ is locally bounded and vanishes at infinity, and also that $\int_{\bbsR^d}h^+(x)^\xi dx<+\infty$.
Notice also that the following condition $({\bf C'}_\xi)$ implies condition $({\bf C}_\xi)$:
\begin{enumerate}
\item[$({\bf C'}_\xi)$\quad] there is some $\gamma>d/\xi $ and $C>0$ such that $ |h(x)|\leq C(|x|^{-\gamma}\wedge 1),\ x\in\bbR^d$.
\end{enumerate}
We will need the following Lemma:
\begin{lemme}\label{lem3} Let $\xi>0$ and $\bar G\in  RV_{-\xi}$.
\begin{enumerate}
\item If  $h$ satisfies condition $({\bf C}_\xi)$, then, for all $\lambda>0$, $\alpha^+(h,\lambda G)=0$ and also $\alpha^+(h,G_\xi)=0$.
\item $\alpha^+(h,G_\xi)=0$ if and only if $\int_{\bbsR^d}h^+(x)^\xi dx<+\infty$.
\end{enumerate}
\end{lemme}
In view of this Lemma, Propositions \ref{prop:stat}, \ref{prop:pr1}, \ref{prop:pr2} and  \ref{prop:pr3} imply:
\begin{corol}
Suppose that $\bar G\in RV_{-\xi}$ and that $h$ is a continuous and satisfies condition $({\bf C}_\xi)$. 
Then for all $\lambda>0$, the random field $M_\lambda\sim {\rm ESN}(h,\lambda G)$  is stationary, $\alpha$-mixing,  almost surely finite and continuous.  
\end{corol}

{\it Proof of Lemma \ref{lem3}:}
We have to check that for any $u>0$, 
$$\int_{\bbsR^d} \bar G(u/h^+(x))dx <+\infty .$$
Let $\delta\in (0,\xi)$ as given by condition $({\bf C}_\xi)$. Since $\bar G\in RV_{-\xi}$, then  there is $C>0$ such that
$\bar G(x)\leq C x^{-(\xi-\delta)}$. Note indeed that the function $x\mapsto x^{\xi-\delta}\bar G(x)$ is bounded on $\bbR_+$ since it is equal to $0$ when $x=0$, is càd-làg, and has limit $0$ as $x\to+\infty$. Then, we have 
$$\bar G(u/h^+(x))\leq C h^+(x)^{\xi-\delta}u^{-(\xi-\delta)}$$ 
and condition $({\bf C}_\xi)$ ensures that these functions are integrable on $\bbR^d$. This shows that $\alpha^+(h,\lambda G)=0$.
The second point is straightforward since 
$$\int_{\bbsR^d} \bar G_\xi(u/h^+(x))dx=u^{-\xi}\int_{\bbsR^d} h^+(x)^\xi dx $$
is finite if and only if $\int_{\bbsR^d} h^+(x)^\xi dx$ is finite. \CQFD 

\begin{remark} If we assume furthermore that the tail function $\bar G$ is such that
$$0<\liminf_{x\to+\infty}x^{-\xi}\bar G(x)\leq \limsup_{x\to+\infty}x^{-\xi}\bar G(x)<+\infty,$$
then condition $({\bf C}_\xi)$ can be replaced by  $\int_{\bbsR^d}h(x)^\xi dx<+\infty$. Lemma \ref{lem3} and Theorems \ref{theo:cv} and \ref{theo:pp} below remain true. The proofs are almost the same and the details will be omitted. 
\end{remark}

\subsection{The large intensity scaling}
We consider the asymptotic behavior of the extremal shot noise $M_\lambda$ as the intensity $\lambda$ goes to infinity.

\begin{theo}\label{theo:cv}
Suppose that $\bar G\in RV_{-\xi}$, that $h$ is continuous and satisfies condition $({\bf C}_\xi)$. Then the following weak convergence of random fields holds in the space $\bbC(\bbR^d,[0,+\infty))$:
$$a_\lambda^{-1}M_\lambda \Longrightarrow {\rm ESN}(h,G_\xi) \quad {\rm as} \ \lambda\to+\infty.$$ 
\end{theo}

\begin{remark}\label{rk:5} Using the terminology in \cite{DM}, we see that under the assumptions of Theorem \ref{theo:cv}, $M_\lambda$ is a regularly varying $\bbC$-valued random field with exponent $\xi>0$.
\end{remark}

{\it Proof of Theorem \ref{theo:cv}:}
We first prove convergence of finite dimensional distributions. From Proposition \ref{prop:pr1}: for $y_1,\ldots,y_k\in\bbR^d$ and $u_1,\ldots,u_k>0$, 
\begin{eqnarray*}
&&\bbP(a_\lambda^{-1}M_\lambda(y_1)\leq u_1,\ldots,a_\lambda^{-1}M_\lambda(y_k)\leq u_k)\\
&=&\exp\left(-\int_{\bbsR^d}\bar G_\lambda\left(\wedge_{1\leq i\leq k}\{u_i/h(y_i-x)\}\right)dx\right)
\end{eqnarray*}
where $\bar G_\lambda(u)=\lambda \bar G(a_\lambda u)$. Equation \eqref{eq:asytail} states that for all $u>0$,  
$\bar G_\lambda(u)\to \bar G_\xi(u)$ as $\lambda+\to\infty$. Hence, as $\lambda\to+\infty$
\begin{eqnarray}
&&\bbP(a_\lambda^{-1}M_\lambda(y_1)\leq u_1,\ldots,a_\lambda^{-1}M_\lambda(y_k)\leq u_k) \nonumber\\
&\to&  \exp\left(-\int_{\bbsR^d}\bar G_\xi\left(\wedge_{1\leq i\leq k}\{u_i/h(y_i-x)\}\right)dx\right), \label{eq:limit}
\end{eqnarray}
provided  Lebesgue's dominated convergence theorem can be applied in order to justify the convergence.
Notice that the right hand side of \eqref{eq:limit} is the cdf of ${\rm ESN}(h,G_\xi)$.
We now justify the convergence (\ref{eq:limit}). Condition $({\bf C}_\xi)$ implies that $h$ is
bounded from above, so that there is $\varepsilon>0$ such that for all $x\in\bbR^d$, 
$$\wedge_{1\leq i\leq k}\{u_i/h(y_i-x)\}>\varepsilon.$$
Then applying Lemma \ref{lem:tail}, there is some $C>0$ such that for large enough $\lambda$
$$\bar G_\lambda(\wedge_{1\leq i\leq k}\{u_i/h(y_i-x)\}) \leq C (\wedge_{1\leq i\leq k}\{u_i/h(y_i-x)\})^{\delta-\xi},\quad x\in\bbR^d.
$$
Condition $({\bf C}_\xi)$ ensures that the right hand side of the above inequality is integrable with respect to $x\in\bbR^d$. Hence equation \eqref{eq:limit} is proved thanks to dominated convergence. This proves the convergence of finite dimensional distributions.

Next, we prove weak convergence in the space $\bbC(\bbR^d,\bbR)$. We prove that for all $u>0$ the weak convergence $(a_\lambda^{-1}M_\lambda)\vee u \Longrightarrow M_\infty\vee u$ holds in $\bbC(\bbR^d,\bbR)$, where $M_\infty\sim{\rm ESN}(h,G_\xi)$. Since $\|(a_\lambda^{-1}M_\lambda)\vee u -a_\lambda^{-1}M_\lambda \|_\infty\leq u$ and $\|M_\infty\vee u -M_\infty \|_\infty\leq u$, this implies the weak convergence $a_\lambda^{-1}M_\lambda \Longrightarrow M_\infty$ (see  \cite{Bill}).

Let $u>0$ be fixed. We have already proved the convergence of the finite dimensional distributions $(a_\lambda^{-1}M_\lambda)\vee u \stackrel{fdd}{\longrightarrow} M_\infty\vee u$. It remains to prove tightness. Condition $({\bf C}_\xi)$ ensures that the function $h$ goes to zero at infinity, and hence $h$ is uniformly continuous on $\bbR^d$. For $\gamma>0$, the  modulus of continuity of $h$ is defined by
$$\omega(h,\gamma)=\sup\{|h(y)-h(x)|;\ |x-y|\leq \gamma\}<+\infty.$$
Let $K=\{z\in\bbR^d; \|z\|\leq 1\}$ be the closed unit ball. Denote by $\Phi_\lambda\sim{\rm PPP}(\lambda dxG(dm))$ the Poisson Point process associated with the extremal shot noise $M_\lambda$. Using Lemma \ref{lem:lem2}, for $y\in K$, we have
\begin{equation}\label{eq:repres2}
(a_\lambda^{-1} M_\lambda(y))\vee u=\sup\{a_\lambda^{-1}m_\phi h(y-x_\phi);\ \phi\in\Phi_{\lambda,K,u}\}\vee u,
\end{equation}
where $\Phi_{\lambda,K,u}=\{\phi\in\Phi_\lambda; \sup_{y\in K}h_\phi(y)\geq a_\lambda u\}$ is almost surely finite. Define
$$T_{\lambda,K,u}=\max\{a_\lambda^{-1}m_\phi; \phi\in \Phi_{\lambda,K,u}\}.$$
The modulus of continuity of $(a_\lambda^{-1}M_\lambda)\vee u$ on $K$ is defined by
$$\omega_K((a_\lambda^{-1}M_\lambda)\vee u,\gamma)=\sup\{|a_\lambda^{-1}M_\lambda(y)\vee u-a_\lambda^{-1}M_\lambda(x)\vee u|;\ x,y\in K,\ |x-y|\leq \gamma\}.$$
Using equation \eqref{eq:repres2} and the definition of $T_{\lambda,K,u}$, the modulus of continuity  satisfies 
$$\omega_K((a_\lambda^{-1}M_\lambda)\vee u,\gamma)\leq T_{\lambda,K,u} \omega(h,\gamma).$$
Hence it is enough to prove that the family $T_{\lambda,K,u}$ is tight. For $v>0$,
\begin{eqnarray*}
\bbP( T_{\lambda,K,u}>v)&=&\bbP(\exists \phi\in \Phi; |a_\lambda^{-1}m_\phi|>v {\rm\ and\ } a_\lambda^{-1} \sup_{y\in K}h_\phi(y)\geq u\}\\
&=& 1-\exp\left(-\int_{\bbsR^d\times (0,+\infty)} 1_{\{a_\lambda^{-1}m>v\}}1_{\{a_\lambda^{-1}mh^+(-x)\geq u\}}\lambda G(dm)dx\right)\\  
&=&1-\exp\left(-\int_{\bbsR^d} \bar G_\lambda(u/h^+(-x)\vee v ) dx\right)\\
&\leq &\int_{\bbsR^d} \bar G_\lambda(u/h^+(-x)\vee v )dx\\
&\leq & C \int_{\bbsR^d} (h^+(-x)^{\xi-\delta}u^{-(\xi-\delta)}\wedge v ^{-(\xi-\delta)} )dx
\end{eqnarray*}
where the last inequality holds for some $C>0$ and $\lambda$ large enough (see Lemma \ref{lem:tail}). This last upper bound is uniform in $\lambda$ and condition $({\bf C}_\xi)$ implies that it goes to zero as $v\to+\infty$. As a consequence,  $T_{\lambda,K,u}$ is tight in $\bbR$, and $(a_\lambda^{-1}M_\lambda)\vee u$ is tight in $\bbC(K,\bbR)$. 
We conclude that $a_\lambda^{-1}M_\lambda\Longrightarrow M_\infty$ in $\bbC(K,\bbR)$.
Using stationarity, the result holds in $\bbC(K+h,\bbR)$ for all $h\in\bbR^d$, and hence in  $\bbC(\bbR^d,\bbR)$.
\CQFD

\subsection{A point process approach}
In this section, we develop a point process framework for the convergence of heavy-tailed extremal shot noise.  Let $\Phi_\lambda\sim {\rm PPP}(\lambda dxG(dm))$ with $\bar G\in RV_{-\xi}$. We show that in a suitable space of functions, the empirical point process
$N_\lambda=\sum_{\phi\in\Phi_\lambda} \delta_{a_\lambda^{-1}h_\phi}$
convergse as $\lambda\to +\infty$, and recover as a by-product Theorem \ref{theo:cv} as well as the convergence of order statistics. For the definition and properties of point processes on general Polish spaces (i.e. complete separable metric spaces), the reader should refer to  \cite{DVJ}.\\ 
\ \\
We start with the presentation of the suitable function spaces adapted from \cite{dH}; note that besides the case of random processes (d=1), this framework covers also the case of random fields ($d\geq 1$). Let $K\subset \bbR^d$ be a compact set and denote by $\bbC^+(K)= \bbC(K,\bbR_+)$ the space of non-negative continuous functions on $K$ endowed with the norm $\|f\|_K=\sup\{|f(y)|;\ y\in K\}$.
Let $\bbC_1^+(K)=\{ f\in \bbC(K);\ f\geq 0 \ \mbox{ and }\ \|f\|_K=1\}$. Thanks to the transformation $f\mapsto (f/\|f\|_K,\|f\|_K)$, we identify $\bbC^+(K)\setminus \{0\}$ to $\bbC_1^+(K)\times (0,+\infty)$.  
We endow $(0,+\infty)$ with the metric $d(u,v)=|1/u-1/v|$, so that its completion is $(0,+\infty]$. Accordingly, we define 
$\bar \bbC^+(K)= \bbC_1^+(K) \times (0,\infty]$ the completion of $\bbC^+(K)\setminus \{0\}$. Note that 
$\bar \bbC^+(K)$ is a Polish metric space and that a bounded subset of $\bar \bbC^+(K)$ is bounded away from zero in the sense that it is included in $\bbC_1^+(K) \times [\varepsilon,\infty]$ for some $\varepsilon>0$.
\ \\
We consider the empirical point process on $\bar \bbC^+(K)$ defined by
$$N_\lambda=\sum_{\phi\in\Phi_\lambda} \delta_{a_\lambda^{-1}h_\phi},\quad \lambda>0.$$
We use a slight abuse of notation here: points $\phi\in\Phi_\lambda$ such that $h_\phi\equiv 0$ on $K$ should be ignored; or equivalently consider the restriction of $N_\lambda$ to $\bar \bbC^+(K)$. However this gives rise to no confusion. 

\begin{theo}\label{theo:pp}
Under the assumptions of Theorem \ref{theo:cv}, the random point process $N_\lambda$ on $\bar\bbC^+(K)$ weakly converges as $\lambda\to+\infty$ to 
$$N_\infty=\sum_{\phi\in\Phi_\infty}\delta_{h_\phi},$$  
with $\Phi_\infty\sim {\rm PPP}(dxG_\xi(dm))$. 
\end{theo}

This result is strongly linked to the so-called POT (Peak Over Threshold) method used by hydrologists and in extreme value theory. Let $f\in\bbC^+(K)$ be a non zero threshold function. For example, Theorem \ref{theo:pp} implies that the number of points  $\phi\in \Phi_\lambda$ such that $h_\phi\geq a_\lambda f$ on $ K$ has a Poisson distribution with mean asymptotically equivalent to
$$\int_{\bbsR^d} \inf_{y\in K}\left( \frac{h(y-x)}{f(y)}\right)^\xi dx .$$
This point process approach is powerful: for instance Theorem \ref{theo:cv} is easily recovered from Theorem \ref{theo:pp}, and new results for order statistics are also easily derived. More precisely, the order statistics are defined by considering the non-increasing reordering of the values $(h_\phi(y);\ \phi\in\Phi_\lambda)$ and are denoted by $M_{\lambda}^{(1)}(y)\geq M_{\lambda}^{(2)}(y)\geq \ldots\geq M_{\lambda}^{(r)}(y))\geq \ldots, r\geq 1$.
Note that the first order statistic coincides with the maximum, i.e. $M_{\lambda}^{(1)}(y)=M_{\lambda}(y)$.

\begin{theo}\label{theo:orderstat}
Under the assumptions of Theorem \ref{theo:cv}, for all $r\geq 1$, the following weak convergence holds in
$\bbC(\bbR^d,[0,+\infty])^r$:
$$(a_\lambda^{-1}M_\lambda^{(i)})_{1\leq i\leq r}\Longrightarrow (M_\infty^{(i)})_{1\leq i\leq r}$$
where $M_\infty^{(i)}$ is the $i$-th order statistic random field associated to 
$(h_\phi;\ \phi\in\Phi_\infty)$ with $\Phi_\infty\sim {\rm PPP}(dxG_\xi(dm))$. 
\end{theo}

Before proceeding to the proof, we recall some notion on convergence of measures on a Polish metric space (see \cite{DVJ}).
A Borel measure $\nu$ on a Polish metric space is boundedly finite if $\nu(A)<+\infty$ for every bounded Borel set $A$.
We say that a sequence of boundedly finite measures $(\nu_k)$ boundedly converges to a boundedly finite measure $\nu$ if $\nu_k(A)\to\nu(A)$ for each bounded Borel set $A$ with $\nu(\partial A)=0$. 

\noindent
{\it Proof of Theorem \ref{theo:pp}:}
Let $\tilde\Phi_\lambda=T_\lambda(\Phi_\lambda)$ be the image of the Poisson Point Process $\Phi_\lambda$ under the transformation $T_\lambda(x,m)=(x,a_\lambda^{-1}m)$.  This is a Poisson point process on $\bbR^d\times (0,+\infty]$ with intensity measure $\nu_\lambda(dx,dm)=dx G_\lambda(dm)$ which is the image measure of $\lambda dx G(dm)$ under $T_\lambda$. Equation \eqref{eq:asytail} implies that $\nu_\lambda$ boundedly converge to $\nu_\infty(dx,dm)=dx G_\xi(dm)$ as $\lambda\to+\infty$ (the metric on $(0,+\infty]$ ensures that $G_\xi$ is boundedly finite). Consequently, the Poisson Point Process $\tilde\Phi_\lambda$ converge to $\Phi_\infty\sim{\rm PPP}(dxG_\xi(dm))$. Let $\cK\subset \bbR^d$ such that $x\in\cK$ if and only if $\|h(.-x)\|_{K}\neq 0$. Since $h$ is continuous and $K$ is compact, $\cK$ is an open subset of $\bbR^d$. The convergence of the restrictions on $\cK\times (0,+\infty]$ also holds:  $\tilde\Phi_\lambda\cap (\cK\times (0,+\infty])$ converges to $\Phi_\infty\cap (\cK\times (0,+\infty])$. 

Next, consider $\Theta:\cK\times (0,\infty] \to \bar\bbC^+(K)$ the application defined by $\Theta(x,m)= mh(\cdot-x)$. The point process $N_\lambda$ (respectively $N_\infty$) on $\bar\bbC^+(K)$ is the image by $\Theta$ of the point process   $\tilde\Phi_\lambda\cap (\cK\times (0,+\infty])$ (resp. $\Phi_\infty\cap (\cK\times (0,+\infty])$) on  $\cK\times (0,+\infty]$. 
The intensity measures of $\tilde \Phi_\lambda$ and $\Phi_\infty$ are $\nu_\lambda$ and $\nu_\infty$ respectively. Then $N_\lambda$ and $N_\infty$ are Poisson Point Processes on $\bar\bbC^+(K)$ with intensity $\nu_\lambda \Theta^{-1}$ and $\nu_\infty \Theta^{-1}$respectively. To prove Theorem \ref{theo:pp}, it is enough to prove that $\nu_\lambda \Theta^{-1}$ and $\nu_\infty \Theta^{-1}$ are boundedly finite measures and that $\nu_\lambda \Theta^{-1}$ boundedly converge to $\nu_\infty \Theta^{-1}$ as $\lambda\to +\infty$.
To that aim, we show that for any bounded set $A\subset \bar\bbC^+(K)$,
\begin{equation}\label{eq:fini}
\nu_\lambda\Theta^{-1}(A)<+\infty, \quad \nu_\infty\Theta^{-1}(A)<+\infty
\end{equation}
and
\begin{equation}\label{eq:conv}
\nu_\lambda\Theta^{-1}(A)\to \nu_\infty\Theta^{-1}(A) \quad \mbox{if}\quad \nu_\infty\Theta^{-1}(\partial A)=0.
\end{equation}
Note that if $A$ is bounded in $\bar\bbC^+(K)$, then $A\subset S_\varepsilon=\bbC_1^+(K)\times (\varepsilon,+\infty]$ for some $\varepsilon>0$. Observe that $\nu_\lambda\Theta^{-1}(S_\varepsilon)<+\infty$ and $\nu\Theta^{-1}(S_\varepsilon)<+\infty$. Indeed:
\begin{eqnarray*}
\nu_\lambda\Theta^{-1}(S_\varepsilon)&=&\nu_\lambda(\{(x,m)\in \bbR^d\times (0,+\infty); \|mh(\cdot-x)\|_K\geq \varepsilon \})\\
&=& \int_{\bbsR^d\times (0,+\infty)} 1_{\{mh_K(-x)\geq \varepsilon\}} dx G_\lambda(dm)\\
&=& \int_{\bbsR^d} \bar G_\lambda( \varepsilon/h_K(x)) dx,
\end{eqnarray*}
and similarly 
$$\nu_\infty\Theta^{-1}(S_\varepsilon)=\int_{\bbR^d} \bar G_\xi( \varepsilon/h_K(x)) dx.$$
Equation \eqref{eq:fini} is hence equivalent to the fact that $\alpha(h_K,G_\lambda)=\alpha(h_K,G_\xi)=0$ and it is enough to check $\alpha^+(h,G)=\alpha^+(h,G_\xi)=0$. From Lemma \ref{lem3}, condition $({\bf C}_\xi)$ implies that $\alpha^+(h,G)$=0. Furthermore, condition $({\bf C}_\xi)$ implies that $\int_{\bbsR^d}h^+(x)^\xi dx<+\infty$ and hence $\alpha^+(h,G_\xi)=0$. This proves equation \eqref{eq:fini}.\\
Let $\nu_{\lambda,\varepsilon}\Theta^{-1}(\cdot)=\nu_{\lambda}\Theta^{-1}(\cdot\cap S_\varepsilon)$ and $\nu_{\varepsilon}\Theta^{-1}(\cdot)=\nu\Theta^{-1}(\cdot\cap S_\varepsilon)$. These are finite measures on $\bar\bbC^+(K)$ and equation \eqref{eq:conv} is equivalent to the weak convergence $\nu_{\lambda,\varepsilon}\Theta^{-1}\Longrightarrow \nu_{\varepsilon}\Theta^{-1}$ on $\bar\bbC^+(K)$. 
Convergence of the finite dimensional distributions is proven as in Theorem \ref{theo:cv} since for $A=\{f\in\bar\bbC^+(K); f(y_1)> u_1,\ldots, f(y_k)> u_k\}$ 
\begin{eqnarray*}
\nu_{\lambda}\Theta^{-1}(A)&=&\int_{\bbsR^d\times (0,+\infty)}1_{\{mh(y_i-x)> u_i;1\leq i\leq n\}}dx G_\lambda(dm)\\
&=& \int_{\bbsR^d}\bar G_\lambda(\vee_{i=1}^n u_i/h(y_i-x))dx\\
&\rightarrow & \int_{\bbsR^d}\bar G_\xi(\vee_{i=1}^n u_i/h(y_i-x))dx\\
&=&\nu_{\infty}\Theta^{-1}(A).
\end{eqnarray*}
The limit is a consequence of Lebesgue's Theorem and Lemma \ref{lem:tail}. It remains to prove tightness. Let $\delta>0$, we prove that for large enough $M$ and $\lambda>1$, we have 
$$\nu_{\lambda,\varepsilon}(\Theta( [-M,M]^d \times [\varepsilon,+\infty]))\geq 1-\delta $$
where $\Theta( [-M,M]^d \times [\varepsilon,+\infty])$ is compact in $\bar\bbC^+(K)$ as the image of a compact set by the countinous application $T$. We have indeed 
\begin{eqnarray*}
\nu_{\lambda,\varepsilon}(\Theta( [-M,M]^d \times [\varepsilon,+\infty]))
&=&1-\int_{\bbsR^d\times (\varepsilon,+\infty)} 1_{\{mh_K(-x)\geq\varepsilon\}}1_{\{|x|>M\}}dx\bar G_\lambda(dm)\\
&=&1-\int_{\{|x|>M\}}\bar G_\lambda(\varepsilon/h_K(x))dx. 
\end{eqnarray*}
Condition $({\bf C}_\xi)$ and Lemma \ref{lem:tail} imply that this last term goes to $1$ as $M\to+\infty$ uniformly in large $\lambda$ (see the proof of Lemma \ref{lem3}). \CQFD
\ \\
{\it Proof of Theorem \ref{theo:orderstat}:}
Let $\varepsilon>0$ and consider 
$$S_\varepsilon=\bbC_1^+(K)\times (\varepsilon,+\infty]\subset \bar\bbC^+(K).$$
Since $S_\varepsilon$ is bounded with $\nu_\infty \Theta^{-1}(\partial S_\varepsilon)=0$, the set $N_\lambda\cap S_\varepsilon$ is a.s. finite and weakly converges to $N_\infty\cap S_\varepsilon$. Let $\cM_p(\bbC(K))$ be the space of finite point measures on $\bbC(K)$. The mapping 
$$\left\{\begin{array}{ccc}\cM_p(\bbC(K))&\to& \bbC(K)\\ \sum_{i=1}^m \delta_{f_i} &\mapsto& \vee_{i=1}^m f_i \vee \varepsilon \end{array}\right.  $$
is continuous. Similarly, for each $r\geq 1$, the following mapping is continuous: 
$$\Psi_r\left\{\begin{array}{ccc}\cM_p(\bbC(K))&\to& \bbC(K)\\ \sum_{i=1}^m \delta_{f_i} &\mapsto& f^{(r)} \vee \varepsilon \end{array}\right.  $$
where $f^{(r)}(y)$ is the $r$-th order statistic in $\{f_i(y); 1\leq i\leq m\}$ and is $0$ if $r>m$.
Then, Theorem \ref{theo:pp} and the continuous mapping Theorem yield the weak convergence on $\bbC(K)$:
$$(a_\lambda^{-1}M_\lambda^{(r)}) \vee \varepsilon = \Psi_r (N_\lambda \cap S_\varepsilon) \Longrightarrow\   \Psi_r(N_\infty \cap S_\varepsilon) = M_\infty^{(r)} \vee \varepsilon.$$
Letting $\varepsilon \to 0$, we have $a_\lambda^{-1}M_\lambda^{(r)}\Longrightarrow M_\infty^{(r)}$ in $\bbC(K)$. The compact $K$ being arbitrary, this proves the convergence of the $r$-th order statistic in $\bbC(\bbR^d,[0,+\infty])$.
In order to consider several order statistics, apply the continuous mapping Theorem to $N_\lambda\cap S_\varepsilon\mapsto (\Psi_i( N_\lambda\cap S_\varepsilon))_{1\leq i\leq r}$.\CQFD
\ \\

\subsection{Supremum of heavy-tailed ESN over large balls}
In this section, we consider the asymptotic behavior of the supremum $\sup_{|y|\leq R} M_\lambda(y)$ as $R\to +\infty$. 

\begin{theo}\label{theo2}
Let $M_\lambda\sim{\rm ESN}(h,\lambda G)$ and suppose that $\bar G\in RV_{-\xi}$ for some $\xi>0$ and  that $h$ satisfies condition $({\bf C'}_\xi)$. Then the following weak convergence holds as $R\to+\infty$,
$$\frac{ 1}{N(R)}\sup_{|y|\leq R} M_{\lambda}(y) \Longrightarrow F_\xi,$$
where $N(R)=\|h\|_\infty c_d^{1/\xi}\lambda^{1/\xi}G^{\leftarrow}(1-R^{-d})$ and $c_d$ denotes the volume of the euclidean unit ball  in $\bbR^d$.
\end{theo}
{\it Proof of Theorem \ref{theo2}:}
Corollary \ref{cor1}, applied to $h_R(x)=\sup\{h(x+z); |z|\leq R\}$, yields
\begin{equation}\label{eq:fbound1}\bbP(\sup_{|y|\leq R} M_{\lambda}(y)\leq u)=\exp\left(-\lambda\int_{\bbsR^d}\bar G\left(u/h_{R}(x)\right)dx\right),\quad u>0.
\end{equation}
Since $h$ satisfies condition $({\bf C'}_\xi)$, $h_R$ satisfies also condition $({\bf C'}_\xi)$ and according to Lemma \ref{lem3}, $\alpha^+(h_R,G)=0$. Furthermore $h$ is bounded and vanishes at infinity, so that there is some $R_0$ such that $\sup_{|y|\leq R_0}h(y)=\|h\|_\infty$. For any $R\geq R_0$, $h_{R}(x)=\|h\|_\infty$ if $|x|\leq R-R_0$ and $h_{R}(x)\leq C((|x|-R)^{-\gamma}\wedge 1)$ if
$|x|\geq R-R_0$. Put $a_{R^d}=G^{\leftarrow}(1-R^{-d})$. We estimate
\begin{eqnarray}
&&\int_{\bbsR^d}\bar G\left(a_{R^d} \|h\|_\infty u/h_{R}(x)\right)dx \nonumber \\
&=& c_d(R-R_0)^d\bar G\left(a_{R^d}u\right)+\int_{|x|>R-R_0}\bar G\left(a_{R^d}\|h\|_\infty u/h_{R}(x)\right)dx.\label{eq:fbound2}
\end{eqnarray}
Equation \eqref{eq:asytail} implies $c_d(R-R_0)^d\bar G(a_{R^d}u)\to c_du^{-\xi}$ as $R\to+\infty$. On the other hand, using that $h_R$ is bounded from above and Lemma \ref{lem:tail}, we get that for large enough $R$,
\begin{equation}\label{eq:fbound3}
\bar G\left(a_{R^d}\|h\|_\infty u/h_{R}(x)\right)\leq C R^{-d} (u/h_R(x))^{-\xi+\delta}. 
\end{equation}
As a consequence, 
\begin{eqnarray}
& &\int_{|x|>R-R_0}\bar G\left(a_{R^d}\|h\|_\infty u/h_{R}(x)\right)dx\nonumber\\
&\leq& C u^{-\xi+\delta} R^{-d} \int_{|x|>R-R_0}  h_R(x)^{(\xi-\delta)}dx\nonumber\\
&\leq& C   u^{-\xi+\delta} R^{-d} \int_{|x|>R-R_0}   (|x|-R)_+^{-\gamma(\xi-\delta)}\wedge 1) dx\nonumber\\
&\leq& C   u^{-\xi+\delta} R^{-d} \int_0^{+\infty} ((r-R_0)_+^{-\gamma(\xi-\delta)}\wedge 1)c_{d-1}(r+R-R_0)^{d-1}dr\label{eq:fbound4}
\end{eqnarray}
The last line is obtained using polar coordinates. Note that the integral in \eqref{eq:fbound4} is finite for $\delta$ small enough since $d<\gamma\xi$ under condition $({\bf C'}_\xi)$. Since the bound \eqref{eq:fbound4} goes to $0$ as $R\to+\infty$,
$$\bbP\left( M_{R,\lambda}\leq a_{R^d} \|h\|_\infty u\right)\to \exp(-\lambda c_du^{-\xi}),\quad {\mbox as} \ R\to+\infty$$
follows from equations \eqref{eq:fbound1}, \eqref{eq:fbound2} and \eqref{eq:fbound3}. This achieves the proof. \CQFD

\section{Properties of the max-stable extremal shot noise}\label{sec:4}
We consider the extremal shot noise $M_\infty$ appearing in Theorem \ref{theo:cv} and investigate its properties.
\begin{prop}\label{prop:pr4}
Let $M_\infty\sim {\rm ESN}(h,G_\xi)$ and suppose $h$ is continuous and such that 
$$\int h^\xi(x) dx=1\quad {\rm  and} \int h^+(x)^\xi dx<+\infty.$$ 
Then  $M_\infty$ is a continuous, stationary, $\alpha$-mixing, max-stable random field with Fréchet margins $F_\xi$.
\end{prop}
The result follows from Propositions \ref{prop:pr2} and \ref{prop:pr3}. The condition $\int h^\xi(x) dx=1$ ensures the normalization to unit Fréchet margins, otherwise a scale parameters appears.\\

A first insight into the dependence structure of the random field $M_\infty$ is given by the extremal coefficient function. 
Recall from \cite{SchTaw} that the extremal coefficient function $\theta(h)$ of the stationary max-stable process $M_\infty$ is defined by the equation
$$\bbP\left[M_\infty(x)<u,M_\infty(x+h)<u\right]=\bbP\left[M_\infty(0)<u\right]^{\theta(h)}.$$
More generally, the extremal coefficient associated with the compact set $K$ is defined by the relation
$$\bbP\left[\sup_{x\in K}M_\infty(x)<u\right]=\bbP\left[M_\infty(0)<u\right]^{\theta(K)}.$$
Furthermore, for each $v\in\bbR^d\setminus\{0\}$, the sequence $(M_\infty(nv))_{n\geq 0}$ is a stationary max-stable sequence with Fréchet marginals $F_\xi$. We denote by $\gamma(v)$ the extremal index of this sequence given by the relation
$$\lim_{n\to+\infty}\bbP\left(\vee_{1\leq i\leq n} M_\infty(nv)\leq n^{1/\xi}u \right)=F_\xi(u)^{\gamma(v)},\quad u>0. $$ 
Note that $\gamma(v)\in[0,1]$ with $\gamma(v)=1$ in the independent case. See \cite{FHR} for a general discussion on the extremal index. 

\begin{prop}\label{prop:pr5}
Under the assumptions of Proposition \ref{prop:pr4}, the extremal coefficient is given by
$$\theta(K)=\int_{\bbsR^d}h_K^\xi(x)dx $$ 
and in particular the extremal coefficient function is
$$\theta(h)=\int_{\bbsR^d}\left(h^\xi(x)\vee h^\xi(x+h)\right) dx.$$
The extremal index in direction $v$ is given by
$$\gamma(v)= \inf_{n\geq 1}\frac{1}{n}\int_{\bbsR^d} \vee_{k=1}^{n} h^\xi(x+kv)dx.$$ 
\end{prop}
{\it Proof of Proposition \ref{prop:pr5}:}
The extremal coefficient $\theta(K)$ is computed as follows: by Corollary \ref{cor1}
\begin{eqnarray*}
\bbP\left[\sup_{x\in K}M_\infty(x)<u\right]&=& \exp\left(-\int_{\bbsR^d}(u/h_K(x))^{-\xi}dx\right)\\
&=& \bbP\left[M_\infty(0)<u\right]^{\theta(K)}
\end{eqnarray*}
with $\theta(K)=\int_{\bbsR^d}h_K(x)^{\xi}dx$. 

In the case when $K=\{0,h\}$, we have $h_K(x)=h(x)\vee h(x+h)$ which yields 
$\theta(h)=\int_{\bbsR^d}\left(h(x)^{\xi}\vee h(x+h)^{\xi}dx\right)$.\\

In order to compute the extremal index $\gamma(v)$, we remark that for $u>0$,
\begin{eqnarray*}
\bbP\left(\vee_{1\leq i\leq n} M_\infty(nv)\leq n^{1/\xi}u \right)&=& \exp\left(-n^{-1}u^{-\xi}\int_{\bbsR^d} \vee_{k=1}^n h^\xi(x+kv)dx  \right)\\
&=& F_\xi(u)^{\gamma_n(v)}
\end{eqnarray*}
with $\gamma_n(v)=n^{-1}\int_{\bbsR^d} \vee_{k=1}^n h^\xi(x+kv)dx$. We remark that the sequence $n\gamma_n(v)$ is  
subadditive because 
$$\vee_{k=1}^{n+m} h^\xi(x+kv)\leq \vee_{k=1}^{n} h^\xi(x+kv) + \vee_{k=1}^{m} h^\xi(x+kv+ nv).$$
Hence $\gamma_n(v)$ converges to $\gamma(v)=\inf_{n\geq 1} \gamma_n(v)$.\CQFD
\\ \ \\
In the sequel, we consider the structure of extremal points associated to the max-stable random field $M_\infty$.
A point $\phi\in\Phi_\infty$ is said to be extremal if there is some $x\in\bbR^d$ such that $M_\infty(x)=h_\phi(x)$. 
The subset of extremal points is denoted by $\tilde\Phi_\infty$ and satisfies
$$M_\infty(y)=\sup\{h_\phi(y); \phi\in\tilde\Phi_\infty\},\quad y\in\bbR^d.$$
This is the smallest subset of $\Phi_\infty$ with this property.

\begin{prop}\label{prop:pr6}
Under the assumptions of Proposition \ref{prop:pr4}, the point process $\tilde\Phi_\infty$ is a  stationary marked point process on $\bbR^d\times (0,\infty)$. Its Campbell-Matthes measure $C$ defined by
$$C(P\times Q)= \bbE\left[\sum_{\phi\in\tilde\Phi_\infty}1_{\{x_\phi\in P\}}1_{\{m_\phi\in Q\}}\right],\quad P\subset\bbR^d, Q\subset (0,+\infty),$$
satisfies $C(dx,dm)=\tilde\lambda dx \nu(dm)$ where $\tilde\lambda$ is the intensity 
$$\tilde\lambda=\int_{0}^\infty \left(1-\bbP(M_\infty\geq mh)\right) G_\xi(dm)$$ 
and $\nu$ is the Palm distribution of the marks
$$\nu(dm)=\tilde\lambda^{-1}\left(1-\bbP(M_\infty\geq mh)\right)G_\xi(dm).$$ 
\end{prop}
{\it Proof of Proposition \ref{prop:pr6}:}
This is an application of Campbell's formula for Poisson Point Processes. Let $\phi=(x_\phi,m_\phi)$ and $\mu(d\phi)=dxG_\xi(dm)$ the intensity measure of the point process $\Phi$. We have
\begin{eqnarray*}
C(P\times Q)&=& \bbE\left[\sum_{\phi\in\tilde\Phi_\infty}1_{\{x_\phi\in P\}}1_{\{m_\phi\in Q\}}\right] \\
&=& \bbE\left[\sum_{\phi\in\Phi_\infty}1_{\{\phi\in P\times Q\}} F(\phi,\Phi_\infty\setminus\{\phi\})\right]
\end{eqnarray*}
with 
$$F(\phi,\Phi)=1_{\{h_\phi \not < \sup_{\psi \in \Phi} h_\psi \}}.$$
Applying Campbell's formula,
\begin{eqnarray*}
C(P\times Q)&=&\bbE\left[ \int 1_{\{\phi\in P\times Q\}}  F(\phi,\Phi_\infty)\mu(d\phi)\right]\\
&=& \int 1_{\{\phi\in P\times Q\}}  \bbP(h_\phi \not < M_\infty)\mu(d\phi)\\
&=& \int 1_{\{x\in P\}}1_{\{m\in Q\}}(1-\bbP( M_\infty>m h(\cdot-x))) dx G_\xi(dm) 
\end{eqnarray*}
where the expectation is taken with respect to $\Phi_\infty$. By stationarity, $\bbP( M_\infty>m h(\cdot-x))$ does not depend on $x$. Hence, we obtain
$$C(dx,dm)= (1-\bbP( M_\infty>m h(\cdot-x))) dx G_\xi(dm)$$
and the intensity measure of the point process and the Palm distribution of the marks are easily deduced.\CQFD

\begin{remark}\label{rk:7} 
In the same way, higher moment measures can  be explicited in terms of  $\bbP(M_\infty>f), f\in\bbC$. For example, if $A_1$ and $A_2$ are disjoint subsets in $\bbR^d\times (0,\infty)$, 
$$C^{(2)}(A_1\times A_2)=\bbE\left[\sum_{(\phi_1,\phi_2)\in\tilde\Phi_\infty^2}1_{\{\phi_1\in A_1\}}1_{\{\phi_2\in A_2\}}\right]$$
satisfies
\begin{eqnarray*}
&&C^{(2)}(A_1\times A_2)\\
&=&\int_{A_1\times A_2}\bbP(h_{\phi_1}\not < M_\infty  \ \mbox{and}\ h_{\phi_2}\not < M_\infty\ )1_{\{h_{\phi_1}\not< h_{\phi_2}\}}1_{\{h_{\phi_2}\not< h_{\phi_1}\}}\mu(d\phi_1)\mu(d\phi_2).
\end{eqnarray*}
\end{remark}
%\begin{remark}\label{rk:8}
%For $K$ a compact subset of $\bbR^d$, we define $\tilde\Phi_\infty^K$ the set of extremal points relatively to $K$ by the property $\phi\in \tilde\Phi_\infty^K$ if and only if there is some $x\in K$ such that $M_\infty(x)=h_\phi(x)$. 
%Lemma \ref{lem:lem2} implies that for any compact $K$, $\tilde\Phi_\infty^K$ is almost surely finite. The point process $\tilde\Phi_\infty^K$ is not stationary but the  results of Proposition \ref{prop:pr6} are easily adapted: for example 
%the intensity measure is given by $f_K(\phi)d\phi$ where $f_K(\phi)=1-\bbP(M_\infty> h_\phi\ \mbox{ on}\ K).$
%When $K=\{z\}$, one can show that $\tilde\Phi_\infty^{\{z\}}$ consists almost surely in a single random point $\phi_z=(x_z,m_z)$ such that $x_{z}$ and $m_{z}h(z-x_{z})$ are independent random variables with distribution $h^\xi(z-x)dx$ and $F_\xi$ respectively.
%\end{remark}

\vspace{1cm}

\appendix
\section{A technical lemma}
\begin{lemme}\label{lem:tail}
Suppose $\bar G\in RV_{-\xi}$ for some $\xi>0$. For $u>0$, define
$$\bar G_\lambda(u)=\lambda \bar G(a_\lambda u),\quad {\rm with}\quad a_\lambda=G^\leftarrow(1-\lambda^{-1}).$$
For all $\delta\in (0,\xi)$ and $\varepsilon>0$, there exist $C>0$ and $\lambda_0>0$ such that
for all $ u\geq\varepsilon$ and $\lambda\geq \lambda_0$,
$$\bar G_\lambda(u)\leq C u^{-(\xi-\delta)}.$$ 
\end{lemme}
{\it Proof:}  Since $\bar G\in RV_{-\xi}$, $a_\lambda=G^\leftarrow(1-\lambda^{-1})\to +\infty$  and 
$c_\lambda=\lambda \bar G(a_\lambda)\to 1$ as $\lambda\to +\infty$. As a consequence, we have
$$\bar G_\lambda(u)=\lambda \bar G(a_\lambda u)= c_\lambda \frac{\bar G(a_\lambda u)}{\bar G(a_\lambda)}.$$ 
Then the Lemma follows from the Potter's bound for the regularly varying function $\bar G\in RV_{-\xi}$ (see \cite{Res} Proposition 0.8 (ii) or \cite{BGT}).\CQFD

\end{document}